\documentclass[11pt]{article}
\usepackage{amsmath,amsfonts,amssymb, amsthm,enumerate,graphicx,xcolor,url}
\usepackage{mathtools}
\usepackage{tikz,authblk}
\usepackage{subfig}

\newtheorem{Theorem} {{\textsf{Theorem}}}[section]
\newtheorem{Proposition}[Theorem]{{\textsf{Proposition}}}
\newtheorem{Corollary}[Theorem]{{\textsf{Corollary}}}

\newtheorem{definition}[Theorem]{{\textsf{Definition}}}
\newtheorem{remark}[Theorem]{{\textsf{Remark}}}

\newtheorem{Lemma}[Theorem]{{\textsf{Lemma}}}

\def\RP2{{\mathbb{R}P}^2}
\def\D{\Delta}
\def\wrt{{\rm with respect to }}

\newcommand{\Star}{\mbox{\upshape st}\,}
\newcommand{\lk}{\mbox{\upshape lk}\,}

\begin{document}

\title{A characterization of $g_2$-minimal normal 3-pseudomanifolds with at most four singularities}
\author{Biplab Basak$^1$, Raju Kumar Gupta and Sourav Sarkar}

\date{}

\maketitle

\vspace{-15mm}
\begin{center}

\noindent {\small Department of Mathematics, Indian Institute of Technology Delhi, New Delhi 110016, India$^{2}$.}

\end{center}

\footnotetext[1]{{\em Corresponding author:}  Biplab Basak}

\footnotetext[2]{{\em E-mail addresses:} \url{biplab@iitd.ac.in} (B.
Basak), \url{Raju.Kumar.Gupta@maths.iitd.ac.in} (R. K. Gupta), \url{Sourav.Sarkar@maths.iitd.ac.in} (S. Sarkar).}

\begin{center}
\date{February 08, 2024}
\end{center}

\hrule

\begin{abstract}
Let $\Delta$ be a $g_2$-minimal normal 3-pseudomanifold. A vertex in $\Delta$ whose link is not a sphere is called a singular vertex. When $\Delta$ contains at most two singular vertices, its combinatorial characterization is known \cite{BasakSwartz}. In this article, we present a combinatorial characterization of such a $\Delta$ when it has three singular vertices, including one $\mathbb{RP}^2$-singularity, or four singular vertices, including two $\mathbb{RP}^2$-singularities. In both cases, we prove that $\Delta$ is obtained from a one-vertex suspension of a surface, and some boundary complexes of $4$-simplices by applying the combinatorial operations of types connected sums, vertex foldings, and edge foldings.
\end{abstract}

\noindent {\small {\em MSC 2020\,:} Primary 57Q05; Secondary 05E45, 51B10, 51E20, 52B70.

}

\medskip

\section{Introduction}
The study of pseudomanifolds, and in particular normal pseudomanifolds, has been one of the central topics in combinatorial topology. Structures for such objects from a combinatorial and topological viewpoint in terms of simplicial and cell complexes are being developed by various researchers throughout the world. A basic and fundamental enumerative invariant of any $n$-dimensional simplicial complex, $\Delta$, is its $f$-vector $(f_{-1},f_0,f_1,\dots, f_n)$, where $f_i$ is the number of $i$-dimensional faces of $\Delta$ for $-1\leq i\leq n$. The empty set is considered to be the only face of dimension $-1$. Inequalities in the lower and upper bound conjectures involving $f_i$ have been used as a very interesting and powerful tool to analyze the interplay between the $f$-vector of a complex and its geometric carrier. In 1970, a breakthrough in this direction was given by Walkup \cite{Walkup}. Inspired by \cite{Walkup}, a new combinatorial invariant $g_2(\D)$ has been defined as  $g_2(\D)=f_1-(d+1)f_0+ {d+2\choose 2}$ for every $d$-dimensional simplicial complex $\D$.  From \cite{Walkup}, it is known that for a triangulated $3$-manifold $\Delta$, $g_2(\Delta)\geq 0$ with the equality occurring when $\Delta$ is a triangulation of a stacked sphere. Barnette \cite{Barnette1, Barnette2} achieved the same results in higher dimensions. In 1987, Kalai \cite{Kalai} used the concept of rigidity and framework on graphs to establish the lower bounds of $f$-vector for normal pseudomanifolds of dimension at least three, where 2-dimensional links were triangulated spheres.

Kalai's findings establish a lower bound on the value of $g_2$ for a normal $d$-pseudomanifold, where 2-dimensional links are spheres, expressed in terms of the same for a link. Subsequently, Fogelsanger's results in [8, Chapter 8] removed the restriction on 2-dimensional links. The combined implications of Kalai and Fogelsanger's results assert that if $d\geq 0$ and $\Delta$ is a normal $d$-pseudomanifold, then for any face $\sigma$ of $\Delta$ with co-dimension 3 or more, the inequality $g_2(\Delta)\geq g_2(\lk (\sigma,\D))$ holds. In \cite{Gromov}, the non-negativity of the invariant $g_2$ has been described in a different way of looking at rigidity, and in \cite{BagchiDatta}, a detailed study on the lower bound theorem has been given in terms of Gromov rigidity. For the case of the  normal 3-pseudomanifold $\Delta$, the result turned out to be $g_2(\Delta)\geq g_2(\lk (v,\D))$ for every vertex $v$ in $\Delta$. If $g_2(\Delta)=g_2(\lk (v,\D))$ for some vertex $v$, then such a complex $\Delta$ is said to be {\em $g_2$-minimal} (or have relatively minimal $g_2$ \cite{BasakSwartz}). Basak and Swartz \cite{BasakSwartz} proved that if $\Delta$ is a $g_2$-minimal normal 3-pseudomanifold with at most two singular vertices, then $\Delta$  is obtained from a one-vertex suspension of a surface and some boundary complexes of $4$-simplices by a sequence of operations of the form vertex foldings and connected sums. In this article, we have worked with the $g_2$-minimal normal $3$-pseudomanifolds with at most four singularities. More precisely, we proved the following:

\begin{Theorem}\label{main theorem}
Let $\Delta$ be a normal $3$-pseudomanifold such that $\Delta$ has $n$ singular vertices among which $n-2$ have $\mathbb{RP}^2$-singularities, where $3\leq n\leq 4$. Then $\Delta$ is obtained from a one-vertex suspension of a surface, and some boundary complexes of $4$-simplices  by a sequence of operations of types connected sums, vertex foldings, and edge foldings.
\end{Theorem}

Several classifications and descriptions of triangulated manifolds and pseudomanifolds have been done based on the $f$-vector and the value of $g_2$. In \cite{BagchiDatta98}, the $d$-pseudomanifold with $d+4$ vertices has been described. The work in \cite{Swartz2008} established that only a finite number of combinatorial manifolds are possible for a given upper bound of $g_2$.
In \cite{Swartz2009, NovikSwartz}, the $f$-vector was the key factor where homology manifolds and specific types of pseudomanifolds were involved. In \cite{NevoNovinsky} and \cite{Zheng}, normal $d$-pseudomanifolds with $g_2=1$ and $g_2=2$ have been characterized, respectively. The concepts of rigidity and the frameworks described in \cite{Kalai, TayWhiteWhiteley} are utilized in each case. In the recent article \cite{BGS}, normal 3-pseudomanifolds with at most two singular vertices are characterized.  Additionally, \cite{BasakGupta} delves into the characterization of normal 3-pseudomanifolds with $g_2 \leq 4$. For a characterization of homology $d$-manifolds with $g_2 \leq 3$, one can refer to \cite{BS}.

\section{Preliminaries}
A {\em simplicial complex} $K$ is a finite collection of simplices in $\mathbb{R}^m$ for some $m\in\mathbb{N}$, such that for every simplex $\sigma\in K$, all of its faces are in $K$, and for any two simplices $\sigma , \tau\in K$, $\sigma\cap \tau$ is a face of both $\sigma$ and $\tau$. We assume that the empty set $\emptyset$ (considered a simplex of dimension $-1$) is a member of every simplicial complex. For a simplicial complex $\Delta$, the notation $|\D|$ will be used to denote the geometric carrier, i.e., $|\D|:= \cup_{\sigma \in \D} \sigma$.  Two simplices, $\sigma = u_0u_1\cdots u_k$ and $\tau = v_0v_1\cdots v_l$ in $\mathbb{R}^n$ for some $n\in \mathbb{N}$, are considered skew if $u_0,\dots,u_k,v_0,\dots,v_l$ are affinely independent. In that case, $u_0\cdots u_kv_0\cdots v_l$ is a $(k + l + 1)$-simplex and is denoted by $\sigma\star\tau$ or $\sigma\tau$. Two simplicial complexes $\Delta_1$ and $\Delta_2$ in some $\mathbb{R}^m$ are called {\em skew} if $\sigma$ and $\tau$ are skew for every simplex $\sigma\in\Delta_1$ and $\tau\in\Delta_2$. The join of two skew simplicial complexes $\D_1$ and $\D_2$ is defined as $\{\sigma \star \tau : \sigma\in \D_1, \tau\in \D_2\}$, and is denoted by $\D_1\star\D_2$. In particular, the join $\sigma\star\D_1$, where $\sigma$ is an $i$-simplex and $\D_1$ is a $j$-dimensional simplicial complex, is defined as the $(i+j+1)$-dimensional simplicial complex $\{\tau :\tau\leq \sigma\}\star\D_1$. If a simplex $\tau$ is a face of another simplex $\sigma$, then we denote it by $\tau\leq \sigma$. If $S$ is a collection of simplices, then the sub-collection of all the simplices of dimension at most one in $S$ is called the {\em graph} of $S$ and is denoted by $G(S)$. The link of a face $\sigma$ in $\D$ is the collection of faces in $\D$ that do not intersect $\sigma$ and whose join with $\sigma$ lies in $\D$ and is denoted by $\lk (\sigma, \D)$. The star of a face $\sigma$ in $\D$ is  $\{\gamma : \gamma\leq \sigma\star\alpha$ and $ \alpha\in \lk (\sigma, \D)\}$ and is denoted by $\Star (\sigma, \D)$. If the underlying simplicial complex is specified, we then simply write $\lk (\sigma)$ and $\Star (\sigma)$ for $\lk (\sigma, \Delta)$ and $\Star (\sigma,\D)$, respectively. A maximal face in $\Delta$ (i.e., a maximal simplex in $\Delta$) is called a {\em facet} of $\Delta$, and if all the facets are of the same dimension, then we say $\Delta$ is a pure simplicial complex. For a simplex $\sigma$, its boundary complex is defined as $\{\tau\,:\, \tau \leq \sigma $ and $\tau \neq \sigma\}$.

 
Several combinatorial tools have been defined to classify simplicial complexes, and taking the suspension of a simplicial complex is one such tool. In \cite{BagchiDatta} and \cite{Bjorner}, another type of suspension has been defined for the class of normal pseudomanifolds by introducing one extra vertex instead of two. Let $\D$ be a normal $(d-1)$-pseudomanifold, and $v\in\D$ be a vertex. Then, for a new symbol $u\notin\D$, consider the  normal $d$-pseudomanifold $\sum_{v,u}\D:=(v\star\{\tau\in\D : v\nleq\tau\})\cup(u\star\D)$. The geometric carrier of $\sum_{v,u}\D$ is the suspension of $|\D|$. The complex $\sum_{v,u}\D$ is called the {\em one-vertex suspension} of $\D$ \wrt the vertex $v$. Note that in $\sum_{v,u}\D$, $\lk u$ and $\lk v$ are simplicially isomorphic to $\D$.
 
%
%
%
 
 A $d$-pseudomanifold is a pure $d$-dimensional simplicial complex $\Delta$ such that every face of dimension $d-1$ is contained in exactly two facets, and for any pair of facets $\sigma$ and $ \tau$ in $\D$, there is a sequence of facets $\sigma_1, \sigma_2,\dots,\sigma_m$ in $\D$ such that $\sigma=\sigma_1, \tau=\sigma_m$, and $\sigma_i\cap\sigma_{i+1}$ is a $(d-1)$-simplex in $\D$ for $1\leq i\leq m-1$. The $d$-pseudomanifolds in which the links of all faces of dimensions up to $d-2$ are connected are called {\em normal}. In particular, a normal 3-pseudomanifold is a strongly connected $3$-dimensional pure simplicial complex where the geometric carrier of the link of each vertex is a closed and connected surface. A vertex in a normal pseudomanifold $\D$ whose link is a triangulated sphere is called a {\em non-singular vertex} of $\D$. If $t$ is a vertex in a normal pseudomanifold $\D$ whose link is not a triangulated sphere, we call $t$ a {\em singular vertex} of $\D$, and if $|\lk (t,\D)|\cong S$, we say $t$ has $S$-singularity.  Two combinatorial tools, vertex folding, and edge folding, are first introduced in \cite{BasakSwartz}, and these two are going to be the main key factors in our result.

 \begin{definition}
Let $\sigma_1$ and $\sigma_2$ be two facets of a simplicial complex $\Delta$ whose intersection is a single vertex $x$. Suppose there is a simplicial bijection $\psi : \sigma_1\to\sigma_2$ such that $\psi (x)=x$ and for all other vertices $y\leq\sigma_1$, the only path of length two or less from $y$ to $\psi (y)$ is $(y,x,\psi (y))$. Then we can obtain a complex $\Delta_{x}^{\psi}$ by identifying all the faces $\rho\leq\sigma_1$ with $\psi (\rho)\leq\sigma_2$ and then removing the identified facet. We call $\Delta_{x}^{\psi}$ a vertex folding of $\Delta$ at $x$, and in this case, $\Delta$ is called a vertex unfolding of $\Delta^{\psi}_{x}$.
 \end{definition}
 \begin{definition}
Let $\sigma_1$ and $\sigma_2$ be two facets of a simplicial complex $\Delta$ whose intersection is a single edge $xy$. Suppose there is a simplicial bijection $\psi : \sigma_1\to\sigma_2$ such that $\psi (x)=x$, $\psi(y)=y$, and for all other vertices $u\leq\sigma_1$, all paths of length two or less from $u$ to $\psi (u)$ pass through either $x$ or $y$. Then we can obtain a complex $\Delta_{xy}^{\psi}$ by identifying all the faces $\rho\leq\sigma_1$ with $\psi (\rho)\leq\sigma_2$ and then removing the identified facet. We call $\Delta_{xy}^{\psi}$ an edge folding of $\Delta$ at $xy$, and in this case, $\Delta$ is called an edge unfolding of $\Delta^{\psi}_{xy}$.
 \end{definition}

If $\Delta$ is a normal $d$-pseudomanifold, then $\Delta^{\psi}_{xy}$ is a pseudomanifold that need not be normal. Moreover,
$$g_2(\Delta^{\psi}_{xy})=g_2(\Delta)+ {d \choose 2},\qquad \text{and}$$
$$\lk (x,\Delta^{\psi}_{xy})=\lk (x,\Delta)^{\psi}_{y}.$$
After applying vertex folding and edge folding to a simplicial complex, the vertices of two facets get identified by the admissible map $\psi$ in the new complex. However, to avoid ambiguity in notations in the new complex, we denote the identified vertices with the same notation as in the original complex.

The definition of a simplicial complex implies that sometimes $\Delta$ may contain all the proper faces of a simplex $\sigma$, but $\sigma \notin \Delta$. In such a case, we say $\sigma$ is a missing face of $\Delta$. If $S$ is a subset of $\Delta$ such that $xy \in S$ but the vertices $x$ and $y$ are not in $S$, then we say that $xy$ is an open edge of $S$ and denote it by $(x, y)$. The presence of a missing tetrahedron and an open edge in a normal 3-pseudomanifold and its subsets provides information regarding the construction of the same.

%

\begin{Proposition}[\cite{BasakSwartz}]\label{Edgefolding}
Let $\Delta$ be a  normal $3$-pseudomanifold. Suppose $pqrs$ is a missing facet in $\Delta$ such that $\partial(qrs)$ and $\partial(prs)$ separate $\lk (p,\Delta)$ and $\lk (q,\Delta)$, respectively, and a small neighborhood of $|\partial(pqr)|$ in $|\lk (s,\Delta)|$ is a M\"{o}bius strip. Then a small neighborhood of $|\partial(pqs)|$ in $|\lk (r,\Delta)|$ is a M\"{o}bius strip. Further, there exists a normal $3$-pseudomanifold $\Delta '$ such that $\Delta$ is obtained from $\Delta '$ by an edge folding at $rs\in\Delta'$, and $pqrs$ is the image of the removed facet.
\end{Proposition}

\begin{Proposition}[\cite{BasakSwartz}]\label{Advanced}
Let $\D$ be a $g_2$-minimal normal $3$-pseudomanifold with at most two singular vertices, $t$ and $s$, such that $g_2(\lk (t,\D))\geq g_2(\lk (s,\D))$. Then $\D$ is obtained from a one-vertex suspension of $\lk (s,\D)$ and some boundary complexes of $4$-simplices by a sequence of operations of the form vertex foldings and connected sums. In particular, if $\D$ has exactly one singular vertex, then  $\D$ is obtained from some boundary complexes of $4$-simplices by a sequence of operations of the form vertex foldings and connected sums.
\end{Proposition}

In Proposition \ref{Advanced}, it is evident that if $\lk (t,\D)\cong \lk (s,\D)\#_{n}H$, where $H$ is either Torus or Klein bottle, then $\D$ is obtained from a one-vertex suspension of $\lk (s,\D)$ and some boundary complexes of 4-simplices by exactly $n$ vertex foldings and finitely many connected sums.

%
%
%
 
 \begin{Proposition} [\cite{BasakSwartz}] \label{lemma:graphlink}
Let $\Delta$ be a $g_2$-minimal normal $3$-pseudomanifold.  Then there exists a normal $3$-pseudomanifold $\bar \Delta$ and a vertex $t \in \bar \Delta$ such that $\Delta$ is obtained from $\bar \Delta$ by (possibly zero) facet subdivisions, and $G(\bar \Delta)=G(\Star(t,\bar \Delta))$.
  \end{Proposition}

\section{Normal $3$-pseudomanifolds with three or four singularities}

Let $\Delta$ be a normal $3$-pseudomanifold and $uv$ be an edge in  $\Delta$. Then  $\lk (v, \D)\setminus  \{\sigma \in \lk (v, \D) : u \leq \sigma\}$ is a triangulation of a surface with boundary. Let $K$ be a triangulation of a surface with a boundary. An edge $e$ is called an {\em interior edge} of $K$ if $e$ is contained in exactly two triangles in $K$, and an edge  $e$ is called a {\em boundary edge} of $K$ if $e$ is contained in exactly one triangle in $K$. A vertex $v$ is called a {\em boundary vertex} of $K$ if there is a boundary edge of $K$ containing $v$. A vertex $u$ is referred to as an {\em interior vertex} of $K$ if there is no boundary edge of $K$ containing $u$.

\begin{Lemma}\label{interior vertex}
Let $\Delta$ be a normal $3$-pseudomanifold such that $G(\Delta)=G(\Star (t,\Delta))$, and let $v$ be a vertex in $\D$ such that $v\neq t$. Then  $\lk (v, \D)\setminus  \{\sigma \in \lk (v, \D) : t \leq \sigma\}$ does not contain any interior vertex. 
	\end{Lemma}
\begin{proof} If  $\lk (v, \D)\setminus  \{\sigma \in \lk (v, \D)\,:\, t \leq \sigma\}$ has an interior vertex, say $x$, then $x\in \lk (t,\D)$, and $vx$ is an edge in $\Delta$, but $vx\not\in \lk (t,\D)$. This leads to a contradiction.
	\end{proof}

Let $v$ be a singular vertex in $\Delta$. If $|\lk (v,\Delta)|\cong\mathbb{RP}^2$, then  $\lk (v, \D)\setminus  \{\sigma \in \lk (v, \D) : t \leq \sigma\}$ is a M\"obius strip, and by  Lemma \ref{interior vertex},   $\lk (v, \D)\setminus  \{\sigma \in \lk (v, \D) : t \leq \sigma\}$ has no interior vertices. Let $M$ be a M\"obius strip that does not have any interior vertices. An edge $xy\in M$ will be called a {\em cut edge} if the edge $xy$ is an interior edge of $M$, and $xy$ separates $M$ into two portions. Then one portion will be a disc, and the other portion will be another M\"obius strip. If $xy$ is an interior edge of $M$, and $xy$ does not separate $M$, then we say $xy$ is a {\em non-cut edge} of  the M\"obius strip $M$.

\begin{Lemma}\label{mstrip1}
Let $\Delta$ be a normal $3$-pseudomanifold such that $G(\Delta)=G(\Star (t,\Delta))$  and $\Delta$ has exactly three singular vertices, including $t$. Let $v$ be a singular vertex in $\Delta$ such that $v\neq t$ and $|\lk (v,\D)|\cong \mathbb{RP}^2$. Then  $\lk (v, \D)\setminus  \{\sigma \in \lk (v, \D) : t \leq \sigma\}$ contains a  non-cut edge $xy$, where $x$ and $y$ are non-singular vertices in $\Delta$. Moreover, the non-cut edge $xy$ creates a missing tetrahedron $xyvt$ in $\D$ such that a small neighborhood of $|\partial(xyt)|$ in $|\lk (v,\D)|$ is a M\"obius strip.
\end{Lemma}
\begin{proof}
From Lemma \ref{interior vertex}, we know that the M\"obius strip  $\lk (v, \D)\setminus  \{\sigma \in \lk (v, \D) : t \leq \sigma\}$ does not contain any interior vertex. If the  M\"obius strip  $\lk (v, \D)\setminus  \{\sigma \in \lk (v, \D) : t \leq \sigma\}$ has a cut edge, then we can cut along the cut edge, and we get another M\"obius strip. Continuing this a finite number of times, we have a  M\"obius strip $M$ where each interior edge will be a non-cut edge. Note that a non-cut edge of $M$ is also a non-cut edge of  $\lk (v, \D)\setminus  \{\sigma \in \lk (v, \D) : t \leq \sigma\}$. Furthermore, any vertex $x\in M$ is on the boundary of $M$ and has a degree of at least 3. Moreover, $x$ is adjacent to two vertices in the boundary of $M$, and for the other vertex $y$, if $xy$ is an edge in $M$, then $xy$ is an interior edge, and hence a non-cut edge of $M$. Let $p$ be the singular vertex other than $t$ and $v$ in $\Delta$. If possible, let $M$ not contain an interior edge $xy$, where $x$ and $y$ are non-singular vertices in $\Delta$. Then, all the interior edges of $M$ are incident to $p$ in $M$. If a non-singular vertex $u$ is not adjacent to $p$ in $M$, then $d(u)\geq 3$ implies that there is a non-singular vertex $x$ such that $ux$ is an interior edge of $M$, leading to a contradiction. Thus, all non-singular vertices of $M$ must be adjacent to $p$ in $M$. Consequently, the M\"obius strip $M$ is the same as $\Star(p, M)$. This is contradictory, as $\Star(p, M)$ is a 2-ball. Therefore, $M$ contains an interior edge $xy$, where $x$ and $y$ are non-singular vertices in $\Delta$. Since each interior edge of $M$ is a non-cut edge of $M$, and a non-cut edge of $M$ is also a non-cut edge of $\lk (v, \D)\setminus  \{\sigma \in \lk (v, \D)\,:\, t \leq \sigma\}$, the result is achieved.

It follows from Lemma \ref{interior vertex} that $x,y \in  \lk (tv, \D)$. Since $xy$ is an interior edge of $\lk (v, \D)\setminus  \{\sigma \in \lk (v, \D)\,:\, t \leq \sigma\}$, $xy\not \in \lk (tv, \D)$. Thus, $xyvt \not \in \D$. On the other hand, $xy\in \lk (v, \D)$, $x,y \in  \lk (tv, \D)$, and $G(\Delta)=G(\Star (t,\Delta))$ imply that $\partial(xyvt) \subset \D$. Therefore,  $xyvt$ is a missing tetrahedron in $\D$.
 Since $xy$ is a non-cut edge of $\lk (v, \D)\setminus  \{\sigma \in \lk (v, \D)\,:\, t \leq \sigma\}$, $\partial(txy)$ does not separate $\lk (v,\D)$. Since $|\lk (v,\D)|\cong \mathbb{RP}^2$, a small neighborhood of  $|\partial(txy)|$ in $|\lk (v,\D)|$ is a M\"obius strip.
\end{proof}
\begin{Lemma}\label{mstrip}
Let $\Delta$ be a normal $3$-pseudomanifold such that $G(\Delta)=G(\Star (t,\Delta))$ and $\Delta$ has exactly four singular vertices, including $t$. Let $v$ be a singular vertex in $\Delta$ such that $v\neq t$ and $|\lk (v,\D)|\cong \mathbb{RP}^2$. Then  $\lk (v, \D)\setminus  \{\sigma \in \lk (v, \D) : t \leq \sigma\}$ contains a  non-cut edge $xy$, where $x$ and $y$ are non-singular vertices in $\Delta$. Moreover, the non-cut edge $xy$ creates a missing tetrahedron $xyvt$ in $\D$ such that a small neighborhood of $|\partial(xyt)|$ in $|\lk (v,\D)|$ is a M\"obius strip.
\end{Lemma}
\begin{proof}
Let $p$ and $q$ be two singular vertices other than $v$ and $t$ in $\D$. From Lemma \ref{interior vertex}, we know that the M\"obius strip  $\lk (v, \D)\setminus  \{\sigma \in \lk (v, \D) : t \leq \sigma\}$ does not contain any interior vertex. If the  M\"obius strip  $\lk (v, \D)\setminus  \{\sigma \in \lk (v, \D) : t \leq \sigma\}$ has a cut edge, then we can cut along the cut edge and get another M\"obius strip. Continuing this a finite number of times, we have a  M\"obius strip $M$ where each interior edge is a non-cut edge. Furthermore, any vertex $x\in M$ is in the boundary of $M$ and has a degree of at least 3. Moreover, $x$ is adjacent to two vertices in the boundary of $M$, and for the other vertex $y$, if $xy$ is an edge in $M$, then $xy$ is an interior edge and hence a non-cut edge of $M$. If $M$ contains exactly one singular vertex, then we achieve our result by using similar arguments as in Lemma \ref{mstrip1}.

\begin{figure}[ht]
		\tikzstyle{ver}=[]
		\tikzstyle{vertex}=[circle, draw, fill=black!50, inner sep=0pt, minimum width=2pt]
		\tikzstyle{edge} = [draw,thick,-]
		\centering
		\begin{tikzpicture}[scale=0.45]
			
			\begin{scope}[shift={(-5,0)}]
				\foreach \x/\y/\z in {0/0/p,-2.5/3/p_1,0/4/p_2,2.7/2.5/p_3,2.5/-2.5/p_4,0/-4/p_5,-2/-3/p_6}{
					\node[vertex] (\z) at (\x,\y){};
				}
				\foreach \x/\y/\z in {0.7/0.1/p,-2.7/3.4/p_1,0/4.4/p_2,2.95/2.95/p_3,2.9/-3.1/p_{n-2},0/-4.5/p_{n-1},-2.4/-3.5/p_{n}}{
					\node[ver] () at (\x,\y){$\z$};
				}

				\foreach \x/\y in {p/p_1,p/p_2,p/p_3,p/p_4,p/p_5,p/p_6,p_1/p_2,p_2/p_3,p_4/p_5,p_5/p_6}{
					\path[edge] (\x) -- (\y);}

					\path[edge, dashed] (p_3) -- (p_4);
				
			\end{scope}
			
			\begin{scope}[shift={(5,0)}]
				\foreach \x/\y/\z in {0/0/q,-2.5/3/q_1,0/4/q_2,2.7/2.5/q_3,2.5/-2.5/q_4,0/-4/q_5,-2/-3/q_6}{
					\node[vertex] (\z) at (\x,\y){};
				}
				\foreach \x/\y/\z in {0.7/0.1/q,-2.7/3.4/q_1,0/4.4/q_2,2.95/2.95/q_3,2.9/-3.1/q_{m-2},0/-4.5/q_{m-1},-2.4/-3.5/q_m}{
					\node[ver] () at (\x,\y){$\z$};
				}

				\foreach \x/\y in {q/q_1,q/q_2,q/q_3,q/q_4,q/q_5,q/q_6,q_1/q_2,q_2/q_3,q_4/q_5,q_5/q_6}{
					\path[edge] (\x) -- (\y);}

					\path[edge, dashed] (q_3) -- (q_4);
				
			\end{scope}
			
		\end{tikzpicture}
	\caption{} \label{fig:1}	
	\end{figure}

Assume that $M$ contains both singular vertices $p$ and $q$. If possible, suppose there is no interior edge $xy$ in $M$, where $x$ and $y$ are non-singular vertices. In such a case, all the interior edges of $M$ are incident to either $p$ or $q$. If a non-singular vertex $u$ is adjacent to neither $p$ nor $q$ in $M$, then $d(u)\geq 3$ implies the existence of a non-singular vertex $x$ such that $ux$ is an interior edge of $M$, leading to a contradiction. Thus, all non-singular vertices of $M$ must be connected to either $p$ or $q$ in $M$. Consequently, the  M\"obius strip $M$ is the same as the union of $\Star(p, M)$ and $\Star(q, M)$. Now we consider three possibilities.

\smallskip

\noindent\textbf{Case 1:} Let the edge $pq\not\in M$. Then, the stars $\Star (p, M)$ and $\Star (q, M)$ take the forms depicted in Figure \ref{fig:1}, where $p_i$ and $q_j$ are non-singular vertices, $1\leq i \leq n$, $1\leq j \leq m$.  If we identify $p_{i}p_{i+1}$ with $q_{j}q_{j+1}$ or $q_{j+1}q_{j}$, then the identified edge will be in the interior of $M$, which is a contradiction. Therefore, $\Star (p, M)$ and $\Star (q, M)$ are disjoint, which contradicts the fact that the M\"obius strip $M$ is the same as the union of $\Star (p, M)$ and $\Star (q, M)$.

\begin{figure}[ht]
		\tikzstyle{ver}=[]
		\tikzstyle{vertex}=[circle, draw, fill=black!50, inner sep=0pt, minimum width=2pt]
		\tikzstyle{edge} = [draw,thick,-]
		\centering
		\begin{tikzpicture}[scale=0.4]
			
			\begin{scope}[shift={(0,0)}]
				\foreach \x/\y/\z in {0/1/p,0/-2/q,-3.2/6/p_1,-0.5/6/p_2,2/5.5/p_3,4/2.5/p_4,4/0.7/p_5,4/-1/x,4/-3/q_4,3/-6/q_3,1.5/-8/q_2,-2/-8/q_1}{
					\node[vertex] (\z) at (\x,\y){};
				}
				
				\foreach \x/\y/\z in {-0.8/1.2/p,-0.5/-2.2/q,-2.7/6.5/p_1,-0.5/6.5/p_2,2/6/p_3,5/2.5/p_{n-1},4.7/0.7/p_{n},4.5/-1/x,4.7/-3/q_m,3.5/-6.2/q_3,1.5/-8.5/q_2,-2/-8.5/q_1,0/-10.5/(a)}{
					\node[ver] () at (\x,\y){$\z$};
				}

				\foreach \x/\y in {p/p_1,p/p_2,p/p_3,p/p_4,p/p_5,p/x,p_1/p_2,p_2/p_3,p_4/p_5,p_5/x,q/x,q/q_1,q/q_2,q/q_3,q/q_4,x/q_4,q_3/q_2,q_2/q_1,p/q}{
					\path[edge] (\x) -- (\y);}

					\path[edge, dashed] (p_3) -- (p_4);
					\path[edge, dashed] (q_3) -- (q_4);
				
			\end{scope}

			\begin{scope}[shift={(14,0)}]
		
				\foreach \x/\y/\z in {0/3/p,0/-3/q,3.3/7/p1',4.4/5.5/p2',4.9/4/p3',5.2/1.5/pn2',5.2/-1.5/y,1.7/-8/q1',3.3/-7/q2',4.4/-5.5/q3',4.9/-3.5/qm2',-3.3/-7/q1,-4.4/-5.5/q2,-4.9/-4/q3,-5.2/-1.5/qm1,-5.2/1.2/x,-1.7/8/p1,-3.3/7/p2,-4.4/5.5/p3,-4.98/3/pn1}{
					\node[vertex] (\z) at (\x,\y){};
				}
				
				\foreach \x/\y/\z in {-0.5/2.3/p,-0.5/-2.2/q,4.1/7.3/p_1',5.1/5.9/p_2',5.6/4.2/p_3',6.2/1.7/p_{n_2}',5.6/-1.7/y,2/-8.4/q_1',3.5/-7.4/q_2',4.7/-5.9/q_3',5.9/-4/q_{m_2}',-4.1/-7.3/q_1,-5.1/-5.9/q_2,-5.6/-4.2/q_3,-6.2/-1.7/q_{m_1},-5.7/1.4/x,-2/8.4/p_1,-3.5/7.4/p_2,-4.7/5.9/p_3,-5.8/3.4/p_{n_1}, 0/-10.5/(b)}{
					\node[ver] () at (\x,\y){$\z$};
				}

				\foreach \x/\y in {p/q,p/p1',p/p2',p/p3',p/y,p/pn2',p/p1,p/p2,p/p3,p/pn1,p/x,q/x,q/y,q/qm1,q/q3,q/q2,q/q1,q/q1',q/q2',q/q3',q/qm2',p1/p2,p2/p3,pn1/x,x/qm1,q3/q2,q2/q1,q1'/q2',q2'/q3',qm2'/y,y/pn2',p3'/p2',p2'/p1'}{
					\path[edge] (\x) -- (\y);}

						\foreach \x/\y in {p3/pn1,q3/qm1,q3'/qm2',pn2'/p3'}{
					\path[edge, dashed] (\x) -- (\y);}
				
			\end{scope}
		\end{tikzpicture}
		\caption{} \label{fig:2}
	\end{figure}

\medskip

\noindent\textbf{Case 2:} Let $pq$ be a boundary edge of $M$. Suppose $\lk (p, M)=P(p_1,\dots,p_n,x,q)$. Then,  no edge of the form $qz$ can be identified with some $p_ip_{i+1}$ or $p_{i+1}p_{i}$. Let $\lk (q, M)=P(q_1,\dots,q_m$, $x,p)$.  Here, $x, p_i$, and $ q_j$ are non-singular vertices,  $1\leq i \leq n$, $1\leq j \leq m$.  If we identify $p_{i}p_{i+1}$ with $q_{j}q_{j+1}$ or $q_{j+1}q_{j}$, then the identified edge will be an interior edge of $M$, which is a contradiction. Therefore, the union of $\Star (p, M)$ and $\Star (q, M)$ is a 2-ball (cf. Figure  \ref{fig:2} $(a)$), contradicting the fact that the M\"obius strip $M$ is the same as the union of $\Star (p, M)$ and $\Star (q, M)$. 
	
\medskip

\noindent\textbf{Case 3:} Let $pq$ be an interior edge of $M$. Let $\lk (p, M)=P(p_1,p_2,\dots,p_{n_1},x,q,y,p_{n_2}', p_{n_2-1}'$, $\dots, p_2',p_1')$. Then, no edge of the form $qz$ can be identified with an edge in the path $P(p_{n_1}, p_{n_1-1},\dots,p_1,p,p_1',p_2',\dots,p_{n_2}')$. Let $\lk (q, M)=P(q_1,q_2,\dots,q_{m_1},x,p,y,q_{m_2}', q_{m_2-1}'$, $\dots,q_2',q_1')$.  Here, $x, y, p_i,  p_j', q_k$, and $q_l'$ are non-singular vertices,  $1\leq i \leq n_1$, $1\leq j \leq n_2$, $1\leq k \leq m_1$, $1\leq l \leq m_2$.  If we identify any edge from the paths $P(p_{n_1}, p_{n_1-1},\dots,p_1)$ or $P(p_1',p_2',\dots,p_{n_2}')$ with an edge from the paths $P(q_{m_1}, q_{m_1-1}$, $\dots,q_1)$ or $P(q_1',q_2',\dots,q_{m_2}')$, then the identified edge will be an interior edge of $M$, leading to a contradiction. Consequently, the union of $\Star (p, M)$ and $\Star (q, M)$ is a 2-ball (cf. Figure  \ref{fig:2} $(b)$), contradicting the fact that the M\"obius strip $M$ is the same as the union of $\Star (p, M)$ and $\Star (q, M)$. 

\medskip

Therefore, $M$ contains an interior edge $xy$, where $x$ and $y$ are non-singular vertices in $\Delta$. Since every interior edge of $M$ is a non-cut edge, and a non-cut edge of $M$ is likewise a non-cut edge of $\lk (v, \D)\setminus  \{\sigma \in \lk (v, \D)\,:\, t \leq \sigma\}$, we have successfully established the result.

Using similar arguments as those in Lemma \ref{mstrip1}, it can be concluded that $xyvt$ is a missing tetrahedron in $\D$. Since $xy$ is a non-cut edge of $\lk (v, \D)\setminus  \{\sigma \in \lk (v, \D)\,:\, t \leq \sigma\}$, the boundary $\partial(txy)$ does not separate $\lk (v,\D)$. Since $|\lk (v,\D)|\cong \mathbb{RP}^2$, a small neighborhood of  $|\partial(txy)|$ in $|\lk (v,\D)|$ is a M\"obius strip.
\end{proof}

\begin{Theorem}\label{Main1}
Let $\Delta$ be a $g_2$-minimal normal $3$-pseudomanifold, and let $t$ be a singular vertex in $\D$ such that $G(\Delta)=G(\Star (t,\Delta))$.
\begin{enumerate}
\item If $\D$ has three singular vertices, including one $\mathbb{RP}^2$-singularity, then $\D$ is obtained by an edge folding from a $g_2$-minimal normal $3$-pseudomanifold with exactly two singularities.
\item If $\D$ has four singular vertices, including two $\mathbb{RP}^2$-singularities, then $\D$ is obtained by two repeated edge foldings on a $g_2$-minimal normal $3$-pseudomanifold with exactly two singularities.
 \end{enumerate}
\end{Theorem}
\begin{proof}

First, we assume that $\D$ has three singular vertices, and let the singular vertices of $\D$ other than $t$ be $s$ and $v_1$, where $|\lk (v_1,\D)|\cong \mathbb{RP}^2$. Since $|\lk (v_1,\D)|\cong \mathbb{RP}^2$, by Lemma \ref{mstrip1}, there is a missing tetrahedron $tv_1xy$, where $x$ and $y$ are non-singular vertices. Moreover, a small neighborhood of $|\partial(txy)|$ in $|\lk (v_1,\D)|$ is a M\"obius strip. According to Proposition \ref{Edgefolding}, there exists a normal $3$-pseudomanifold $\D'$ such that $\D=(\D')_{tv_{1}}^{\psi}$. Since $\D$ is $g_2$-minimal \wrt $t$, $\D'$ is also $g_2$-minimal \wrt $t$. On the other hand, $g_2(\lk (v_1,\D'))=g_2(\lk (v_1,\D))-3=0$, and hence $v_1$ is a non-singular vertex in $\D'$. Since $\lk (s,\D)$ is isomorphic to $\lk (s,\D')$, $\D'$ is a $g_2$-minimal normal pseudomanifold with two singular vertices $t$ and $s$. Moreover, $g_2(\D')=g_2 (\lk (t,\D'))$, and $\D$ is obtained from $\D'$ by an edge folding along the edge $tv_1$, where $v_1$ is a non-singular vertex in $\D'$.    

Now, let $\D$ have four singular vertices, including two $\mathbb{RP}^2$-singularities. Let $v_1$ and $v_2$ be the vertices in $\D$ such that $|\lk (v_1,\D)|\cong |\lk (v_2,\D)| \cong\mathbb{RP}^2$. Since $|\lk (v_2,\D)|\cong \mathbb{RP}^2$, by Lemma \ref{mstrip}, there is a missing tetrahedron $tv_2wz$, where $w$ and $z$ are non-singular vertices. Moreover, a small neighborhood of $|\partial(twz)|$ in $|\lk (v_2,\D)|$ is a M\"obius strip. Therefore, by Proposition \ref{Edgefolding}, $\D$ is formed via an edge folding along the edge $tv_2$ from a normal 3-pseudomanifold $\D'$. In $\D'$, $v_2$ is a non-singular vertex because $g_2(\lk (v_2,\D'))=g_2(\lk (v_2,\D))-3=0$. Since the singular vertices other than $t$ and $v_2$ did not participate in the process of edge folding, their link in $\D'$ will be isomorphic to the links in $\D$. On the other hand,  $\D'$ is also a $g_2$-minimal normal pseudomanifold, with $g_2(\D')=g_2(\lk (t,\D'))$. Thus, $\D'$ is a $g_2$ minimal normal pseudomanifold with three singular vertices, and $v_1\in\D'$ is a vertex such that $|\lk (v_1,\D')|\cong \mathbb{RP}^2$. By the same argument as in the last paragraph, $\D'$ is obtained from a $g_2$-minimal normal 3-pseudomanifold $\D''$ with exactly two singular vertices by an edge folding along the edge $tv_1$, where $v_1$ is a non-singular vertex in $\D''$. Therefore, $\D$ is obtained from $\D''$ by two edge foldings along the edges $tv_1$ and $tv_2$, where $v_1$ and $v_2$ are non-singular vertices in $\D''$.
\end{proof}

Let $\D$ be a $g_2$-minimal normal $3$-pseudomanifold with three singular vertices, $t, s$, and $v$, where $|\lk (v,\D)|\cong\mathbb{RP}^2$. Also, let $g_2(\D)=g_2(\lk (t,\D))$. It follows from Proposition \ref{lemma:graphlink} that $\D=\bar{ \D}\#\D_1\#\cdots\#\D_n$, where each $\D_i$ is the boundary complex of a $4$-simplex, and $G(\bar {\Delta})=G(\Star (t,\bar \Delta))$. Note that $\bar \D$ contains a copy of each of the vertices $t, s$, and $v$, with the links having the same Betti numbers as the link of those vertices in $\D$. Let $b_1$ and $b_2$ be the first Betti numbers of $\lk (t,\D)$ and $\lk (s,\D)$, respectively. By Theorem \ref{Main1}, $\bar \D$ is obtained from a $g_2$-minimal normal 3-pseudomanifold $\D'$ by an edge folding along an edge of the form $tv$, where $v$ is a non-singular vertex in $\D'$. 
Further, $\D'$ is $g_2$-minimal \wrt the singular vertex $t\in\D'$, and the first Betti number of the link of $s$ in $\D'$ is unchanged. Notice from the definition that an edge folding along an edge $tx$ reduces the first Betti number of the corresponding vertices by 1 in the resulting complex, and it follows that $b_1-b_2\geq 1$. On the other hand, if $\D$ has four singular vertices, including two $\mathbb{RP}^2$-singularities, and $G(\D)=G(\Star (t,\D))$, then $\D$ is obtained from a $g_2$-minimal normal pseudomanifold by two repeated applications of edge folding. Therefore, by Proposition \ref{lemma:graphlink} and Theorem \ref{Main1}, we have the following:
 
\begin{Corollary}
Let $\D$ be a $g_2$-minimal normal $3$-pseudomanifold with four singular vertices $t, s, v_1,$ and $v_2$, where $|\lk (v_1,\D)|\cong|\lk (v_2,\D)|\cong\mathbb{RP}^2$ and $\D$ is $g_2$-minimal with respect to the vertex $t$. If $b_1$ and $b_2$ are the first Betti numbers of $\lk (t,\D)$ and $\lk (s,\D)$, respectively, then $b_1-b_2\geq 2$.
\end{Corollary}
Let $\D$ be a $g_2$-minimal normal 3-pseudomanifold with three singular vertices, including one $\mathbb{RP}^2$-singularity, or four singular vertices, including two $\mathbb{RP}^2$-singularities. Let $s$ and $t$ be the other singular vertices in $\D$, where $g_2(\D)=g_2(\lk (t,\D))$. Then in the first case, $\lk (t,\D)\cong \lk (s,\D)\#_{n}\mathbb{T}^2\#\mathbb{RP}^2$, and for the second case, $\lk (t,\D)\cong \lk (s,\D)\#_{n}\mathbb{T}^2\#\mathbb{RP}^2\#\mathbb{RP}^2$, where $\mathbb{T}^2$ denotes a triangulated torus. Thus, using Proposition \ref{lemma:graphlink}, Theorem \ref{Main1}, and Proposition \ref{Advanced}, we have the following corollaries:

\begin{Corollary}\label{corollary:1}
Let $\D$ be a $g_2$-minimal normal $3$-pseudomanifold with three singular vertices, including one $\mathbb{RP}^2$-singularity. Let $t$ be the singular vertex in $\D$ such that $g_2(\D)=g_2(\lk (t,\D))$. Then, $\D$ is obtained from a one-vertex suspension of $\lk (s,\D)$ and some boundary complexes of $4$-simplices by $n$ vertex foldings, one edge folding, and finitely many connected sums, where $\lk (t,\D)\cong \lk (s,\D)\#_{n}\mathbb{T}^2\#\mathbb{RP}^2$ for a singular vertex $s$.
\end{Corollary}

\begin{Corollary} \label{corollary:2}
Let $\D$ be a $g_2$-minimal normal $3$-pseudomanifold with four singular vertices, including two $\mathbb{RP}^2$-singularities. Let $t$ be the singular vertex in $\D$ such that $g_2(\D)=g_2(\lk (t,\D))$. Then $\D$ is obtained from a one-vertex suspension of $\lk (s,\D)$ and some boundary complexes of $4$-simplices by $n$ vertex foldings, two edge foldings, and finitely many connected sums, where $\lk (t,\D)\cong \lk (s,\D)\#_{n}\mathbb{T}^2\#\mathbb{RP}^2\#\mathbb{RP}^2$  for a singular vertex $s$.  
\end{Corollary}

\bigskip

\noindent {\em Proof of Theorem} \ref{main theorem}. 
The proof follows from  Corollaries \ref{corollary:1} and \ref{corollary:2}.
 \hfill $\Box$

\bigskip

\begin{remark}
{\rm
The proof of Theorem \ref{main theorem} depends on Lemmas \ref{mstrip1} and \ref{mstrip}, which give the existence of a  non-cut edge $xy$ in  $\lk (v, \D)\setminus  \{\sigma \in \lk (v, \D) : t \leq \sigma\}$, where $v$ is a singular vertex and $x, y$ are non-singular vertices in $\Delta$. If $\Delta$ has more than four singular vertices, then such a non-cut edge may not exist. Therefore, if $\Delta$ has more than four singular vertices, then such a characterization for $\Delta$ may not be possible. In fact, we have an example of a $g_2$-minimal normal 3-pseudomanifold   $\Delta$ with five singular vertices, among which one has $\mathbb{T}^2$-singularity and the remaining four have $\mathbb{RP}^2$-singularities (see \cite{Akhmejanov}), where $\Delta$ cannot be obtained from a normal 3-pseudomanifold   $\Delta'$ by an edge folding.
}
\end{remark}
\bigskip

  \noindent {\bf Acknowledgement:} The authors would like to thank the anonymous referees for their many useful comments and suggestions.  The first author is supported by the Science and Engineering Research Board (CRG/2021/000859). The second author is supported by CSIR (India). The third author is supported by the Prime Minister's Research Fellows (PMRF) Scheme.

\medskip


\end{document}